\documentclass[12pt,a4paper]{article} %[a-shadowing_cor 13.04 - 19.08.22 - 01.01.23]   Updated
\usepackage{amsfonts,amssymb} 
\usepackage[cp1251]{inputenc}

\usepackage[margin=25mm, a4paper]{geometry}  % journal quality

\usepackage{tikz}  

\def\n{\noindent}  \def\?#1{}
\def\IZ{{\mathbb{Z}}}  \def\IR{{\mathbb{R}}}
 \def\cA{{\cal A}}   
      \def\cS{{\cal S}}  
   \def\Tor{{\mathbb{T}}}   \def\la{\lambda}  \def\cN{{\cal N}}

   \def\phi{\varphi}  
\def\ep{\varepsilon}  \def\t{\tilde} \def\bdelete#1{}
    \def\v#1{\vec{#1}}
\def\mod1{\,({\rm mod\ } 1)\,}

\def\beq#1#2{\begin{equation} \label{#1} #2 \end{equation}}
\def\bea#1{\begin{eqnarray*} #1 \end{eqnarray*}} \def\a{\!\!\!&\!\!\!\!&}
\def\toas#1{\stackrel{#1}{\longrightarrow}}
\def\function#1{\left\{\!\!\!\begin{array}{ll} #1 \end{array} \right.}
\def\proof{\smallskip \noindent {\bf Proof. \ }}       %start of proof
\def\blanksquare{\,\,\,$\sqcup\!\!\!\!\sqcap$}         %blank  square
\def\qed{\hfill\blanksquare\linebreak\smallskip\par}   %end of proof
\def\thname{Theorem}  \def\lmname{Lemma}    \def\prname{Proposition}
\def\dfname{Definition}  \def\crname{Corollary}  \def\rmname{Remark}
\def\exname{Example}  

\newtheorem{theorem}{\thname}[section]   %Numbering: Theorem--Other section
\newtheorem{lemma}{\lmname}[section]     %{lemma}[theorem]{Lemma}   subsection
\newtheorem{proposition}[lemma]{\prname} %lemma
\newtheorem{example}{\exname}[section]

\newtheorem{dftn}{\dfname}[section]
\newenvironment{definition}{\begin{dftn}\rm}{\end{dftn}} %section
\def\bdef#1{\begin{definition} #1 \end{definition}}
\newtheorem{rmrk}[lemma]{\rmname}
\newenvironment{remark}{\begin{rmrk}\rm}{\end{rmrk}}     %lemma

\makeatletter\def\fps@figure{htbp}\makeatother %figure pos: tbp - standard

\begin{document}

\title{Average shadowing revisited}
\author{Michael Blank\thanks{
        Institute for Information Transmission Problems RAS
        (Kharkevich Institute);}
        \thanks{National Research University ``Higher School of Economics'';
        e-mail: blank@iitp.ru}
       }
\date{January 2, 2023} %\today}
\maketitle

\begin{abstract} We propose a novel unifying approach to study the shadowing property 
for a broad class of dynamical systems (in particular, discontinuous and non-invertible) 
under a variety of perturbations. In distinction to known constructions, our approach is local: 
it is based on the gluing property which takes into account the shadowing under a single 
(not necessarily small) perturbation. \end{abstract}

{\small\n
}
%http://aps.ecnu.edu.cn/UserFiles/File/MSC2020-Mathematica l% 37Bxx Topological dynamics
% 37B65 Approximate trajectories, pseudotrajectories, shadowing
% 37B05 Dynamical systems involving transformations and group actions with special properties (minimality, distality, proximality, expansivity, etc.)
% 37B10 Symbolic dynamics
% 37C50 Approximate trajectories (pseudotrajectories, shadowing, etc.) in smooth dynamics
%%%%%%%%%%%%%%%%%%%%%%%%%%%%%%%%

\section{Introduction}
When modeling a time-evolving process, we obtain its approximate realizations.
This proximity is due to several reasons. 
First, we never know exactly the description of the process itself, 
and second, the presence of various kinds of errors from purely random 
to rounding errors when implemented on a computer are inevitable. 
The question of the adequacy of the simulation results is primarily 
associated with the presence of a real trajectory of the process under 
study in the vicinity of the obtained realization over the longest possible 
time interval. 
This question is especially nontrivial in the case of a chaotic system, 
since for such systems close trajectories diverge very quickly 
(often exponentially fast).

At the level of correspondence between the individual trajectories of a hyperbolic system and the 
pseudo-trajectories\footnote{Approximate trajectories of a system under small 
     perturbations, already considered by G. Birkhoff \cite[(1927)]{Bi} for a completely different purpose.}, 
this problem was first posed by D.V. Anosov \cite[(1967-70)]{An,An2} as a key step in analyzing 
the structural stability of diffeomorphisms.
A similar but much less intuitive approach, called ``specification'', was proposed in the same 
setting by R. Bowen \cite[(1975)]{Bo}. 
Informally, both approaches guarantee that errors do not accumulate during the modeling process: 
in systems with the shadowing property, each approximate trajectory can be uniformly traced 
by the true trajectory over an arbitrary long period of time. Naturally, this is of great
importance in chaotic systems, where even an arbitrary small error in the initial position leads to 
(exponentially in time) a large divergence of trajectories.

Further development (see main results, generalizations and numerous references in 
monographes dedicated to this subject \cite{Pi, PS} and the textbook \cite{KB}) 
demonstrated deep connections between the shadowing property and various ergodic 
characteristics of dynamical systems. 
In particular, it was shown that for diffeomorphisms the shadowing property implies 
the uniform hyperbolicity. 
To some extent, this restricts the theory of uniform shadowing to an important but very special 
class of hyperbolic dynamical systems.

The concept of average shadowing introduced in \cite{Bl88} about 30 years ago gave a possibility 
to extend significantly the range of perturbations under consideration in the theory of shadowing, 
in particular to be able to deal with perturbations which are small only on average but not 
uniformly. However, the original idea under this concept was twofold: (i) to extend the range 
of perturbations and (ii) to be able to deal with non-hyperbolic systems. 
While the first objective was largely achieved (see discussion of unresolved issues below), 
the second objective was not resolved: the proof was given only for smooth hyperbolic systems.
Nevertheless this concept generated a number of subsequent works where various versions 
of shadowing similar in idea to the average shadowing were introduced (see, for example, 
\cite{KKO,KO,LS,Sa,Sa2,WOC}) and the connections between them were studied in detail.

The most notorious in the variety of obstacles in the analysis of the shadowing property 
is that one needs to take into account an infinite number of independent perturbations of 
the original system. This makes the problem highly nonlocal. It is therefore very desirable 
to reduce the shadowing problem to the situation with a single perturbation, albeit with 
tighter control of the approximation accuracy. 

To realize this idea in our recent paper \cite{Bl22} we developed a fundamentally new 
``gluing'' construction, consisting in the effective approximation of a pair of consecutive 
segments of true trajectories. See exact definitions and details of the construction in Section~\ref{s:pre}. 

In \cite{Bl22} using this construction we were able to study systems under perturbations 
being small in various senses (see Section~\ref{s:pre} for details). Moreover, using it we  
studied some combinations of types of perturbations and types of shadowing.\footnote{Previously, 
    a separate method  was developed to analyze each specific combination of a perturbation 
    and a type of shadowing.}
Still the most interesting Gaussian perturbations were out of reach (it was expected that 
more sophisticated estimates will allow to achieve this objective). In this paper we indeed 
overcome this difficulty.

The paper is organized as follows. In Section~\ref{s:pre} we give general definitions 
related to the shadowing property and introduce the key tool of our analysis -- the gluing 
property. In Section~\ref{s:main} we formulate and prove the main result -- Theorem~\ref{t:main}, 
which deduces various versions of shadowing from the gluing property. 
The remaining part of the paper is devoted to the analysis of the applicability 
of the gluing property for various classes of discrete time dynamical systems, starting with 
hyperbolic diffeomorphisms (Section~\ref{s:diff}), for non-uniformly hyperbolic endomorphisms 
with singularities (Section~\ref{s:endo}), and finally for a special class of multivalued maps 
induced by symbolic dynamics (Section~\ref{s:symb}). 

Among other things, the examples of the systems under study demonstrate the difference between 
strong (\ref{e:glu}) and weak (\ref{e:glu-w}) versions of the gluing property, which shows 
that even with uniformly small perturbations it is possible that only the average shadowing 
takes place (but not the uniform one). Similarly, it is shown that there are systems, 
belonging to the class $\cS(R,A)$, but not to $\cS(A',A)$ 
(see definition of $\cS(\cdot,\cdot)$ in Section~\ref{s:pre}).

\section{Setting and main result}\label{s:pre}
We restrict ourselves to discrete time dynamical systems, leaving the extension of our 
approach to continuous time systems (flows) for future research. A discrete time dynamical 
system is completely defined by a non-necessarily invertible map $T:X \to X$ from a metric 
space $(X,\rho)$ into itself. 

\bdef{A {\em trajectory} of the map $T$ starting at a point $x\in X$ is a sequence 
of points $\v{x}:=\{\dots,x_{-2},x_{-1},x_0,x_1,x_2,\dots\}\subset X$, 
for which $x_0=x$ and $Tx_i=x_{i+1}$ 
for all available indices $i$. The part of $\v{x}$ corresponding to non-negative indices 
is called the {\em forward (semi-)trajectory}, while the part corresponding to non-positive indices 
is called the {\em backward (semi-)trajectory}.}

Observe that although the forward trajectory is always uniquely determined by $x=x_0$ and infinite, 
the backward trajectory might be finite (if its ``last'' point has no preimages)\footnote{In this case we 
   are speaking only about available indices $i$ in the definition.} 
and for a given $x=x_0$ there might be arbitrary many admissible backward trajectories. 

\begin{remark}
Despite the introduction of the backward trajectory of a non-invertible dynamical system 
looks somewhat unusual, we inevitably have to go back and through in time when 
constructing the true trajectory approximating the trajectory of the perturbed system. 
Therefore it is more convenient  to define bi-infinite trajectories from the very beginning.
\end{remark}

\bdef{A {\em pseudo-trajectory} of the map $T$ is a sequence 
of points $\v{y}:=\{\dots,y_{-2},y_{-1},y_0,y_1,y_2,\dots\}\subset X$, for which the sequence of 
distances $\{\rho(Ty_i,y_{i+1})\}$ for all available indices $i$ satisfies a certain ``smallness'' condition. 
The parts corresponding to non-negative or non-positive indices are referred as forward or 
backward pseudo-trajectories.}

Introduce the set of ``moments of perturbations'':
$$\cN(\v{y}):=\{t_i:~ \gamma_{t_i}:=\rho(Ty_{t_i}, y_{t_{i+1}})>0, ~i\in \IZ\}$$ 
ordered with respect to their values, i.e. $t_i<t_{i+1}~\forall i$. We refer to the 
amplitudes of perturbations $\gamma_{t_i}$ as {\em gaps} between consecutive segments 
of true trajectories.

\bdef{For a given $\ep>0$ we say that a pseudo-trajectory $\v{y}$ is of 
\begin{itemize}
\item[(U)] {\em uniform} type, if $\rho(Ty_i,y_{i+1})\le\ep$ for all available indices $i$.
\item[(A)] {\em small on average (strong)} type, if $\exists N$ such that
    $\frac1{2n+1}\sum\limits_{i=-n}^n\rho(Ty_i,y_{i+1})\le\ep~~\forall n\ge N$. 
\item[(A')] {\em small on average (weak)} type, if 
    $\limsup\limits_{n\to\infty}\frac1{2n+1}\sum\limits_{i=-n}^n\rho(Ty_i,y_{i+1})\le\ep$. 
\item[(R)] {\em rare perturbations} type, if the upper density of the set $\cN(\v{y})$ does not exceed $\ep$. 
    Namely, $\limsup\limits_{n\to\infty}\frac1{2n+1} \#(\cN(\v{y})\cap [-n,n]) \le\ep$. 
\end{itemize}}
If the backward pseudo-trajectory is finite, only positive indices $i$ are 
taken into account, which leads to one-sided sums $\frac1{n+1}\sum\limits_{i=0}^n\rho(Ty_i,y_{i+1})$ 
instead of two-sided ones.

The U-type pseudo-trajectory is the classical one, introduced by G.~Birkhoff \cite{Bi} and 
D.V.~Anosov \cite{An}. The A-type was proposed by M.~Blank \cite{Bl88} in order to take 
care about perturbations small only on average. It is unnecessarily strong in the sense that 
the corresponding inequality holds for all sufficiently large $n$. In particular, for true Gaussian 
perturbations, the probability of this event is zero. The reason of this is pure technical: 
this was necessary for the techniques applied in \cite{Bl88} for the analysis of shadowing. 
Thus, the weak version (A'), despite being more natural, is completely new. 
The R-type was introduced as an intermediate version in our recent publication \cite{Bl22}, 
and it allows to consider large but rare perturbations. 

Clearly $U\subset A \subset A'$ and $R \subset A'$, but $R\setminus A\ne\emptyset$. 
Despite $R \subset A'$ we will demonstrate that their separate analysis is worth 
doing since in some situations the R-type perturbations, but not of general A'-type, are shadowed.

To simplify notation we will speak about $\ep$-pseudo-trajectories, when the corresponding 
property is satisfied with the accuracy $\ep$. 

\bigskip

The idea of {\em shadowing} in the theory of dynamical systems boils down to the question,  
is it possible to approximate the pseudo-trajectories of a given dynamical system with true trajectories? 
Naturally, the answer depends on the type of approximation.

\bdef{We say that a true trajectory $\v{x}$ {\em shadows} a pseudo-trajectory $\v{y}$ 
with accuracy $\delta$ (notation $\delta$-shadows): 
\begin{itemize}
\item[(U)] {\em uniformly}, if $\rho(x_i,y_i)\le\delta$ for all available indices $i$.
\item[(A)] {\em on average}, if 
      $\limsup\limits_{n\to\infty} \frac1{2n+1}\sum\limits_{i=-n}^n\rho(x_i,y_i)\le\delta$. 
\end{itemize}}
If the backward pseudo-trajectory is finite, only positive indices $i$ are 
taken into account, which leads to the one-sided sum $\frac1{n+1}\sum\limits_{i=0}^n\rho(x_i,y_i)$.

The U-type shadowing was originally proposed by D.V.~Anosov \cite{An}, 
while the A-type was introduced\footnote{The reason is that pseudo-trajectories with large 
       perturbations cannot be uniformly shadowed.} by M.~Blank \cite{Bl22}.
Naturally, the types of pseudo-trajectories and the types of shadowing may be paired 
in an arbitrary way.

\bdef{We say that a DS $(T,X,\rho)$ satisfies the {\em $(\alpha+\beta)$-shadowing property} 
(notation $T\in \cS(\alpha,\beta)$) 
with $\alpha\in\{U,A,A',R\},~ \beta\in\{U,A\}$ if $\forall\delta>0~\exists\ep>0$ such that 
each $\ep$-pseudo-trajectory of $\alpha$-type can be shadowed in the $\beta$ sense with the 
corresponding accuracy $\delta$.}

For example, $\cS(U,U)$ stands for the classical situation of the uniform shadowing of uniformly 
perturbed systems, while $\cS(A',A)$ corresponds to the average shadowing in the case 
of small on average in the weak sense perturbations.

One of the most interesting open questions related to the shadowing problem is to find out 
under what conditions on the map does the presence of a certain type of shadowing for 
each pseudo-trajectory with a single perturbation implies one or another type 
of shadowing property for the system? 
The reason for this question is that the case of a single perturbation is much simpler, 
and therefore the idea of getting information about other types of perturbations 
from this fact is quite attractive.
The answer is known (although very partially, see \cite{Ach}) only in the case of 
U-shadowing of the so-called positively expansive\footnote{Roughly speaking,  
    this means that if two forward trajectories are uniformly close enough to each 
    other, then they coincide. In particular, this property is satisfied for expanding maps.} 
dynamical systems, if additionally one assumes that the single perturbation does 
not exceed $0<\ep\ll1$. 

In order to give the answer to this question we introduce the following property.

\bdef{We say that a trajectory $\v{z}$ {\em glues}  together semi-trajectories $\v{x}, \v{y}$ 
with accuracy rate $\phi:\IZ\to\IR_+$ {\em strongly} if 
\beq{e:glu}{
     \rho(x_k, z_k)\le\phi(k)\rho(x_0,y_0) ~~\forall k<0, \quad  
     \rho(y_k, z_k)\le\phi(k)\rho(x_0,y_0) ~~\forall k\ge0 }
and {\em weakly} if 
\beq{e:glu-w}{
     \rho(x_k, z_k)\le\phi(k) ~~\forall k<0, \quad  
     \rho(y_k, z_k)\le\phi(k) ~~\forall k\ge0 .}  }
In other words $\v{z}$ approximates both the backward part of $\v{x}$ and the 
forward part of $\v{y}$ with accuracy controlled by the rate function $\phi$, and 
in the strong version the accuracy additionally depends multiplicatively on the 
distance between the ``end-points'' of the glued segments of trajectories.

Without loss of generality, we assume that the functions $\phi(|k|)$ and $\phi(-|k|)$ 
are monotonic. Indeed, replacing a general $\phi$ by its monotone envelope 
$$\t\phi(k):=\function{\sup_{i\le k}\phi(i) &\mbox{if } k<0 \\  
                               \sup_{i\ge k}\phi(i) &\mbox{if } k\ge0} ,$$ 
we get the result. 

\bdef{We say that the DS $(T,X,\rho)$ satisfies the (strong/weak) {\em gluing property} 
with the rate-function $\phi:\IZ\to\IR$ (notation $T\in G_{s/w}(\phi)$) if for any pair of 
trajectories $\v{x}, \v{y}$ there is a trajectory $\v{z}$, which glues them at time $t=0$ 
with accuracy $\phi$ in the strong/weak sense.}

\begin{remark}
If $T\in G_{s/w}(\phi)$, then $\forall \tau\in\IZ$ and for any pair of trajectories $\v{x}, \v{y}$ 
there exists a trajectory $\v{z}$, which glues them (strongly/weakly) at time $t=\tau$ with accuracy $\phi$. 
\end{remark}
Indeed, for a given $\tau$ consider a pair of trajectories $\v{x'}, \v{y'}$ obtained from 
$\v{x}, \v{y}$ by the time shift by $\tau$, namely $x'_i:=x_{i+\tau}, ~y'_i:=y_{i+\tau}, ~\forall i$. 
Then since $\v{x'}, \v{y'}$ may be glued together at time $t=0$ with accuracy $\phi$,  
we deduce the same property for $\v{x}, \v{y}$ by the time $t=\tau$. \qed

\begin{remark} 
Additional assumptions about the rate-function $\phi$ are necessary in order to obtain 
meaningful applications of the gluing property. In what follows, we will only assume  
summability of this function: 
\beq{e:sum}{ \Phi:=\sum_k\phi(k)<\infty .}
\end{remark}

Our main result is the following statement.

\begin{theorem}\label{t:main}
Let $T:X \to X$ be a map from a metric space $(X,\rho)$ into itself, and 
let $T\in G_s(\phi)$ with $\Phi:=\sum_k\phi(k)<\infty$. Then
(a) $T\in \cS(U,U)$, ~~(b) $T\in \cS(A',A)$.
\end{theorem}

\begin{remark}\label{r:bounded}
\begin{enumerate}
\item The part (a) has been proven under the same assumption in \cite{Bl22}, 
         as well, as that $T\in \cS(U,U)$ if $T\in G_w(\phi)$ and the perturbations are 
         uniformly bounded.
        Therefore we will prove only the (b) part here. 
        On the other hand, we will demonstrate that (a) follows from the proof of (b) 
        if one assumes that additionally the perturbations are bounded by $0<\ep\ll1$.
\item The proof of the situation $T\in \cS(R,A)$ in \cite{Bl22} was based on a very crude 
        technical assumption that all the perturbations are equal to some constant $D$. 
        This is not sufficient to get the claim (b), which we do in the present paper.
\item For $\cS(U,U)$ it is enough to check the gluing property for $\v{x}, \v{y}$ with 
        $\rho(x_0,y_0)\le\ep_0\ll1$,
%\item For $\cS(A',A)$ it is enough to check the gluing property for $\v{x}, \v{y}$ with 
%        $\rho(x_0,y_0)\le De^\Phi$,
\item $\cS(U,A)\setminus \cS(U,U) \ne\emptyset$. 
\end{enumerate} 
\end{remark}

\section{Proof of Theorem~\ref{t:main}}\label{s:main}

\proof We will prove a stronger ``linear'' version of the shadowing, 
namely that there is a constant $K=K(\phi)<\infty$, such that for each $\ep>0$ small enough 
there is a true trajectory average approximating with accuracy $\delta\le K\ep$ each 
$\ep$-pseudo-trajectory of A'-type. 

The scheme of the proof goes as follows. For a pseudo-trajectory $\v{y}:=\{y_i\}_{i\in\IZ}$ 
we consider in detail the set of moments of perturbations 
$$\cN(\v{y}):=\{t_i:~ \gamma_{t_i}:=\rho(Ty_{t_i}, y_{t_{i+1}})>0, ~i\in \IZ\} .$$ 
Between the moments of time $t_k$ there are no perturbations and hence $\v{y}$ can 
be divided into segments of true trajectories. Thanks to the $G_s(\phi)$ property each pair of 
consecutive segments of true trajectories can be ``glued'' together by a true trajectory with the 
controlled accuracy. 

Without loss of generality, we will assume that perturbations occur at every moment of time
and therefore $t_i=i~~\forall i\in\IZ$. 

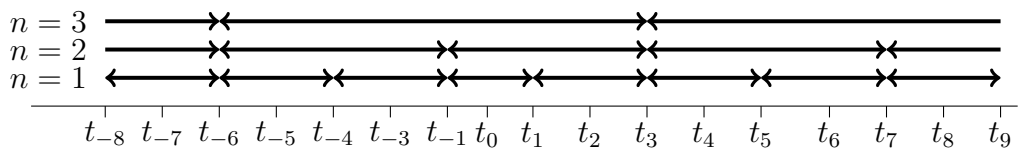
\begin{figure}\begin{center}
\begin{tikzpicture}[scale=0.75]
      \draw [-](-.8,0) to (16.5,0);
      %\draw (-.5,0) to (-.5,-.2); \node at (-.5,-.5){$t_{-9}$};
      \draw (.5,0) to (0.5,-.2); \node at (.5,-.5){$t_{-8}$};
      \draw (1.5,0) to (1.5,-.2); \node at (1.5,-.5){$t_{-7}$};
      \draw (2.5,0) to (2.5,-.2); \node at (2.5,-.5){$t_{-6}$};
      \draw (3.5,0) to (3.5,-.2); \node at (3.5,-.5){$t_{-5}$};
      \draw (4.5,0) to (4.5,-.2); \node at (4.5,-.5){$t_{-4}$};
      \draw (5.5,0) to (5.5,-.2); \node at (5.5,-.5){$t_{-3}$};
      \draw (6.5,0) to (6.5,-.2); \node at (6.5,-.5){$t_{-1}$};
      \draw (7.2,0) to (7.2,-.2); \node at (7.2,-.5){$t_{0}$};
      \draw (8.,0) to (8.0,-.2); \node at (8.0,-.5){$t_{1}$};
      \draw (9.,0) to (9.0,-.2); \node at (9.0,-.5){$t_{2}$};
      \draw (10.,0) to (10.0,-.2); \node at (10.0,-.5){$t_{3}$};
      \draw (11.,0) to (11.0,-.2); \node at (11.0,-.5){$t_{4}$};
      \draw (12.,0) to (12.0,-.2); \node at (12.0,-.5){$t_{5}$};
      \draw (13.2,0) to (13.2,-.2); \node at (13.2,-.5){$t_{6}$};
      \draw (14.2,0) to (14.2,-.2); \node at (14.2,-.5){$t_{7}$};
      \draw (15.2,0) to (15.2,-.2); \node at (15.2,-.5){$t_{8}$};
      \draw (16.2,0) to (16.2,-.2); \node at (16.2,-.5){$t_{9}$};
      \node at (-.5,.5){$n=1$}; \node at (-.5,1.0){$n=2$}; \node at (-.5,1.5){$n=3$};
      %\draw [thick,<->] (6.5,.5) to (8.0,.5); %\node at (1.5,-.5){$t_{-1}$};
      \draw [line width=1.5pt,,<->] (6.5,.5) to (8.0,.5); 
      \draw [line width=1.5pt,,<->] (8.0,.5) to (10.0,.5);
      \draw [line width=1.5pt,,<->] (10.0,.5) to (12.0,.5);
      \draw [line width=1.5pt,,<->] (12.0,.5) to (14.2,.5);
      \draw [line width=1.5pt,,<->] (14.2,.5) to (16.2,.5);
      \draw [line width=1.5pt,,<->] (4.5,.5) to (6.5,.5); 
      \draw [line width=1.5pt,,<->] (2.5,.5) to (4.5,.5); 
      \draw [line width=1.5pt,,<->] (.5,.5) to (2.5,.5); 
      \draw [line width=1.5pt,,<->] (6.5,1.0) to (10.0,1.0); 
      \draw [line width=1.5pt,,<->] (10.0,1.0) to (14.2,1.0); 
      \draw [line width=1.5pt,,<->] (2.5,1.0) to (6.5,1.0); 
      \draw [line width=1.5pt,,->] (.5,1.0) to (2.5,1.0); 
      \draw [line width=1.5pt,,<-] (14.2,1.0) to (16.2,1.0); 
      \draw [line width=1.5pt,,<->] (2.5,1.5) to (10.0,1.5); 
      \draw [line width=1.5pt,,->] (.5,1.5) to (2.5,1.5); 
      \draw [line width=1.5pt,,<-] (10.0,1.5) to (16.2,1.5); 
\end{tikzpicture}\end{center}
\caption{Order of the parallel gluing.}\label{f:sgluing} 
\end{figure}

In our construction (see Fig.~\ref{f:sgluing}) we first glue together pairs of segments 
around the moments of perturbations $t_i$ with even indices: $i_{\pm2k}$, obtaining longer 
segments of the true ``gluing'' trajectories. At each next step, we simultaneously ``glue'' 
together consecutive pairs of already obtained segments. Consequently, at each step of the 
construction, we get a new pseudo-trajectory, consisting of half the number of segments of 
true trajectories (i.e., only half of the perturbation moments remain) with exponentially 
increasing lengths, but with larger gaps between them (compared to the original gaps).
In the limit, we obtain an approximation of the entire pseudo-trajectory. 

This is the scheme of the parallel gluing construction introduced in our previous 
paper \cite{Bl22}. Here we apply the same scheme with an additional advantage 
that on each step of the construction all the lengths of gluing segments coincide 
and those lengths grow exponentially with the step number. 

Another possibility is to use sequential gluing, which can be described as follows. 
Starting from a certain segment of the trajectory (say, between the time moments 
from $t_0$ to $t_1$), we glue it first with the right neighbor, then with the left 
one (or vice versa). Therefore at each step of the construction a new segment 
of the trajectory is glued to the already approximated ones.
In fact, the construction used in \cite{Bl88} to prove the average shadowing property 
for Anosov systems in the terminology given above, is precisely sequential gluing. 
The advantage of the sequential gluing construction is that the corresponding 
calculations are much simpler, but on closer examination it turns out that 
for their application it is necessary to make much stronger assumptions on 
the rate function $\phi$, in particular, that $\phi(\pm1)<1$. 
Even for uniformly hyperbolic systems, this can only be done for so-called 
Lyapunov metric $\rho$, and not for the general one. 
In most of the examples discussed in Sections~\ref{s:diff}--\ref{s:symb}  
the value $\phi(\pm1)$ turns out to be quite large. 

To estimate approximation errors, we find the accuracy of gluing a pair of 
segments of true trajectories: $v_{-N^-}, v_{-N^-+1},\dots,v_{-1}$ and $v_0,v_1,\dots,v_{N^+}$. 
By the property $G_s(\phi)$ there exists a trajectory $\v{z}\subset X$ such that 

$$ \rho(v_k, z_k)\le\phi(k)\rho(Tv_{-1},v_0) ~~\forall k\in \{-N^-,\dots,N^+\}$$
Therefore 
$$ \sum_{k=-N^-}^{N^+}\rho(z_k, v_k) \le \rho(Tv_{-1},v_0) \sum_{k}\phi(k) 
    = \Phi\cdot \rho(Tv_{-1},v_0) .$$
There are several important points here: 
\begin{enumerate}
\item The approximation accuracy depends only on the gap $\rho(Tv_{-1},v_0)$ between the 
        ``end-points'' of the segments of trajectories glued together.
\item After gluing a pair of trajectory segments, the gaps between the end-points in the 
        next step of the procedure may become larger than the original gaps in $\v{y}$. 
\item Each moment of perturbation $t_i$ is taken into account only once during the gluing process.
\item There is an advantage of the (A'+A) case compared with the (R+A) case: all segments of the 
        true trajectories in the process of gluing has the same length (in distinction to the (R+A) case). 
        As we will see, this helps a lot in the analysis.
\end{enumerate}

In \cite{Bl22}, in proving a similar result for the case (R+A), we have used the simplifying assumption 
that all perturbations are uniformly bounded by some constant $D<\infty$. Under this 
assumption, it has been proven that, although the gaps may increase during the gluing process 
they remain uniformly bounded. On the final stage, this fact, together with the assumption 
that perturbations occur very rarely, allow to estimate the average approximation error.

In the present setting, the perturbations are neither uniformly bounded nor sparse. 
Therefore, a more sophisticated approach is needed. Namely, we first show that the average 
values of gaps during the gluing procedure cannot exceed $K\ep$ for some finite  constant $K\ne K(\ep)$. 
Using this estimate, we prove that there exists a finite constant $C$ such that for any $0<\delta\ll1$ 
small enough and $\ep:=C\delta$ for each $\ep$-pseudo-trajectory $\v{y}$ of A' type there is 
a true trajectory $\v{z}$, which on average approximates $\v{y}$ with accuracy $\delta$.
 
Now we are ready to proceed. On the $n$-th step of the process of gluing we have a bi-infinite collection 
of gaps $\{\gamma_{t_i}^{(n)}\}_{i\in\IZ}$. 
By the assumption that the perturbations are small on average, we get 
$$ \limsup_{k\to\infty}\frac1{2k+1}\sum_{i=-k}^k\gamma_{t_i}^{(0)} \le \ep .$$
Our aim is to show that $\exists C\ne C(\ep)$ such that 
$$ \limsup_{k\to\infty}\frac1{2k+1}\sum_{i=-k}^k\gamma_{t_i}^{(n)} \le C\ep \quad\forall n.$$
Using this inequality we can apply the machinery developed for the analysis of the (R+A) case in \cite{Bl22}.

We start with a recursive upper estimate for the gaps:
\beq{e:rec-est}{
  \gamma_{t_i}^{(n+1)} \le \gamma_{t_i}^{(n)} 
                                   + \phi_-^{(n)}\gamma_{t_{i-1}}^{(n)}  + \phi_+^{(n)} \gamma_{t_{i+1}}^{(n)} ,}
where $\phi_\pm^{(n)}=\phi(\pm2^n)$. Indeed, the lengths of the glued segments 
of trajectories on the $n$-th step of the procedure is equal to $2^n$ and 
$\phi_-^{(n)}\gamma_{t_{i-1}}^{(n)}$ is the upper estimate for the approximation error coming 
from the left, while $\phi_+^{(n)} \gamma_{t_{i+1}}^{(n)}$ is the upper estimate for the 
approximation error coming from the right. 

Rewriting (\ref{e:rec-est}) as follows:
$$ \gamma_{t_i}^{(n+1)} \le (\phi_-^{(n)} + \phi_+^{(n)}) \gamma_{t_i}^{(n)} 
             + \left(  (1 - \phi_-^{(n)} - \phi_+^{(n)}) \gamma_{t_i}^{(n)} 
                       + \phi_-^{(n)}\gamma_{t_{i-1}}^{(n)}  + \phi_+^{(n)} \gamma_{t_{i+1}}^{(n)}  \right) ,$$
we see that in the 1st term the factor $\phi_-^{(n)} + \phi_+^{(n)}$ vanishes with $n$, 
while the 2nd term corresponds to the averaging operator of type $v_i \to (1-a-b)v_i + av_{i-1} + bv_{i+1}$, 
the recursive application of which flattens a sequence $\{v_i\}$ to a constant. 

To make this reasoning accurate, we need some calculations. Without loss of generality, 
we assume that the function $\phi$ is even (i.e. $\phi(-k)=\phi(k)~~\forall k$). 
Indeed, replacing a general $\phi$ by $\t\phi(k):=\max(\phi(-k), \phi(k))~~\forall k$, 
we get the result. 

Denote $R_k^{(n)}:= \sum_{i=-k}^k\gamma_{t_i}^{(n)}$. Then using (\ref{e:rec-est}) we get 
%$$ R_k^{(n+1)} = \sum_{i=-k}^k\gamma_{t_i}^{(n+1)} 
%                        \le \sum_{i=-k}^k\gamma_{t_i}^{(n)} 
%                          + \sum_{i=-k}^k \phi_-^{(n)}\gamma_{t_{i-1}}^{(n)} 
%                          + \sum_{i=-k}^k  \phi_+^{(n)} \gamma_{t_{i+1}}^{(n)}   $$
%$$ = R_k^{(n)} + \phi_-^{(n)} \left(\gamma_{t_{-k-1}}^{(n)}  + R_k^{(n)} - \gamma_{t_{k+1}}^{(n)} \right) 
%                       + \phi_+^{(n)} \left(-\gamma_{t_{-k-1}}^{(n)}  + R_k^{(n)} + \gamma_{t_{k+1}}^{(n)}  \right) $$ 
%$$ = (1 + \phi_-^{(n)} + \phi_+^{(n)}) R_k^{(n)} 
%          + (\phi_-^{(n)} - \phi_+^{(n)})(\gamma_{t_{-k-1}}^{(n)}  -  \gamma_{t_{k+1}}^{(n)}) 
%     = (1 + \phi_-^{(n)} + \phi_+^{(n)}) R_k^{(n)}.  $$

\bea{R_k^{(n+1)} \a= \sum_{i=-k}^k\gamma_{t_i}^{(n+1)} 
                             \le \sum_{i=-k}^k\gamma_{t_i}^{(n)} 
                                 + \sum_{i=-k}^k \phi_-^{(n)}\gamma_{t_{i-1}}^{(n)} 
                                 + \sum_{i=-k}^k  \phi_+^{(n)} \gamma_{t_{i+1}}^{(n)} \\ 
                            \a= R_k^{(n)} + \phi_-^{(n)} \left(\gamma_{t_{-k-1}}^{(n)} + R_k^{(n)} 
                                                                              - \gamma_{t_{k+1}}^{(n)} \right) \\
                           \a\qquad\quad+  \phi_+^{(n)} \left(-\gamma_{t_{-k-1}}^{(n)}  + R_k^{(n)} 
                                                            +  \gamma_{t_{k+1}}^{(n)}  \right) \\
                           \a= (1 + \phi_-^{(n)} + \phi_+^{(n)}) R_k^{(n)} 
                                  + (\phi_-^{(n)} - \phi_+^{(n)})(\gamma_{t_{-k-1}}^{(n)}  
                                                                                -  \gamma_{t_{k+1}}^{(n)}) \\
                           \a=  (1 + \phi_-^{(n)} + \phi_+^{(n)}) R_k^{(n)}
                                  \qquad({\rm since}~\phi_-^{(n)} = \phi_+^{(n)}) \\
                           \a\le\dots \le R_k^{(0)} \prod_{i=0}^n  (1 + \phi_-^{(i)} + \phi_+^{(i)})  .}

\begin{lemma}\label{l:inf-prod} For any sequence of nonnegative real numbers $\{b_k\}_{k\ge1}$ we have 
$$ \limsup_{n\to\infty}\prod_{k\ge1}^n(1+b_k) \le e^{\limsup\limits_{n\to\infty}\sum_{k=1}^n b_k} .$$
\end{lemma}
\proof Denote $B_n:=\prod_{k\ge1}^n(1+b_k), ~S_n:=\sum_{k=1}^n b_k$. 
We proceed by induction on $n$. For $n=1$, the question boils down to 
$1+v \le e^v~\forall v\in\IR$. 
$$ \frac{d}{dv} (1+v) = \function{1 >  e^v = \frac{d}{dv}e^v &\mbox{if } v<0 \\  
                                                 1 =  e^v = \frac{d}{dv}e^v &\mbox{if } v=0 \\  
                                                 1 <  e^v = \frac{d}{dv}e^v &\mbox{if } v>0 .}$$
Therefore the graph of $e^v$ lies above the straight line $1+v$ with the only 
tangent point at the origin. 

This implies that $B_1\le e^{S_1}$. Assume that  $B_n\le e^{S_n}$ for some $n\in\IZ_+$ 
and prove the same inequality for $n+1$. We get 
$$ B_{n+1}=(1+b_{n+1})B_n \le (1+b_{n+1}) e^{S_n} \le e^{b_{n+1}}e^{S_n} = e^{S_{n+1}} .$$
Lemma is proven. \qed

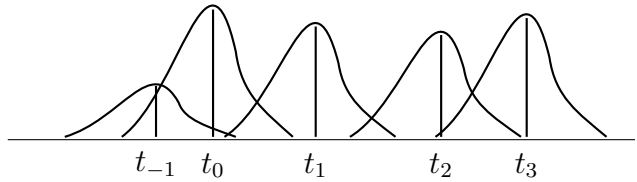
\begin{figure}\begin{center}
\def\G#1#2#3#4{\draw [thick] (1.6+#1,0*#3+#2) to (1.6+#1,1.5*#3+#2); \node at (1.6+#1,-.5+#2){#4};
      \draw [thick] (0+#1,0*#3+#2) .. controls (1+#1,.5*#3+#2) and (1.5+#1,2.5*#3+#2) .. (2+#1,1*#3+#2) 
                                 .. controls (2.1+#1,.5*#3+#2) and (2.5+#1,.3*#3+#2) .. (3+#1,0*#3+#2);}
\begin{tikzpicture}[scale=0.75]
      \draw [-](0,0) to (11,0); 
      \G{1}{.05}{.6}{$t_{-1}$}; \G{2}{.05}{1.5}{$t_{0}$}; \G{3.8}{.05}{1.3}{$t_{1}$};
      \G{6}{.05}{1.2}{$t_{2}$}; \G{7.5}{.05}{1.4}{$t_{3}$};
\end{tikzpicture}\end{center}
\caption{Contributions to the upper bound of the gluing error.}\label{f:error} 
\end{figure}

Now we are ready to continue the proof of Theorem. 
Contributions to the upper bound of the gluing error are coming from two different sources: 
the estimates of the gaps (changing during the steps of the parallel gluing construction) 
and the summation over error contributions from the gluing of pairs of consecutive segments 
of true trajectories (See Fig~\ref{f:error}).

By means of Lemma~\ref{l:inf-prod} we obtain an upper bound for the normalized sum of gaps.
Setting $b_n:=\phi_-^{(n)} + \phi_+^{(n)}$, by Lemma~\ref{l:inf-prod} we get %
\beq{e:est-a}{ R_k^{(n+1)} \le \prod_{i=1}^n(1+b_i) R_k^{(0)} \le e^{\sum_{i=0}^nb_i}R_k^{(0)} 
                        \le e^\Phi R_k^{(0)} .}

For the $n$-th approximating pseudo-trajectory $\v{z}^{(n)}$ denote  
$$ Q_k^{(n)} := \frac1{2k+1} \sum_{t=-k}^k \rho(z_t^{(n)},y_t) .$$
Then %
\bea{Q_k^{(n)} \a\le  \frac1{2k+1} \sum_{t=-k}^k \sum_i \gamma_{-t+i}^{(n)} \phi(i) \\
          \a=  \sum_i \phi(i) \cdot \left( \frac1{2k+1} \sum_{t=-k}^k \gamma_{-t+i}^{(n)}\right)\\
          \a= \sum_i \phi(i) \cdot R_k^{(n)}(t), }
where $R_k^{(n)}(t):=\sum_{i=-k}^k\gamma_{-t+t_i}^{(n)}$.

Therefore, using (\ref{e:est-a}) we obtain an upper bound %
\beq{e:fin-a}{ \limsup_{k\to\infty} Q_k^{(n)}  \le \ep \Phi e^\Phi ,}%
which does not depend on the step number $n$. 

Now for each $k>0$ one estimates the distance between pseudo-trajectories 
$\v{z}^{(n)}$ and $\v{z}^{(n+k)}$ as follows. Recall that the pseudo-trajectory $\v{z}^{(n)}$ 
consists of pieces of true trajectories of length $2^n$. Therefore 
$$ \rho(z^{(n)}_t,z^{(n+k)}_t)  \le \sum_{|j|\ge \tau(t,n)} \phi(j) \gamma_{j}^{(n)} ,$$
where $\tau(t,n)$ for a given $t$ grows to infinity as $2^n$. 

Since the function $\phi$ is summable, this implies that 
$$ \sum_{|j|\ge \tau(t,n)} \phi(j) \toas{n\to\infty} 0 .$$
Therefore, under the additional assumption that the phase space $X$ is compact 
(and hence $\sup_{n,j}\gamma_{j}^{(n)}<\infty$), 
for any given $t$ the sequence $\{z^{(n)}_t\}_n$ 
is fundamental and converges as $n\to\infty$ to the limit $z_t$, where $\{z_t\}$ 
is the true trajectory of our system.

In the non-compact case only weaker average convergence to the limit trajectory $\v{z}$ 
is available: %
$$ \limsup_{m\to\infty}\frac1{2m+1} \sum_{t=-m}^m \rho(z^{(n)}_t,z^{(n+k)}_t) 
     \le \ep e^\Phi \sum_{|j|\ge \tau(t,n)} \phi(j) \toas{n\to\infty} 0 .$$

Since the estimate (\ref{e:fin-a}) is uniform on $n$ we may use it for $\v{z}$ instead 
of $\v{z}^{(n)}$, getting  
$$ \limsup_{k\to\infty}\frac1{2k+1} \sum_{t=-k}^k \rho(z_t,y_t) \le \ep\Phi e^\Phi .$$

Theorem is proven. \qed

\n{\bf Proof of Remark \ref{r:bounded}}. (1) Observe that if the perturbations are uniformly 
bounded by $\ep\ll1$, then from the proof of Theorem~\ref{t:main} on each step of the 
gluing procedure we get that all gaps are bounded from above by $\ep\Phi e^\Phi$, while 
the uniform approximation error is estimated from above by the same constant. 

(3) Follows from the observation above that all gaps are bounded from above by $\ep\Phi e^\Phi$. 

(4) See Proposition \ref{p:non-unif} item (3). \qed

\section{The  gluing property for diffeomorphisms}\label{s:diff}
The gluing property may be explained in terms similar to those which are actively used in the theory 
of smooth hyperbolic dynamical systems. Denote by $\v{x}^-:=\{\dots,x_{-2},x_{-1},x_0=x\}$ and 
$\v{x}^+:=\{x=x_0,x_1,x_x,\dots\}$ be backward and forward semi-trajectories of the point $x\in X$, 
and consider the sets:
$$ W^-(\v{x}^-):=\{\v{z}\subset X:~~\rho(x_k,z_k)\toas{k\to-\infty}0\} .$$ 
$$ W^+(\v{x}^+):=\{\v{z}\subset X:~~\rho(x_k,z_k)\toas{k\to\infty}0\} .$$ 
Then the gluing property means that for each pair of semi-trajectories $\v{x}^-$ and $\v{y}^+$ 
the sets $W^+(\v{x}^+)$ and $W^-(\v{y}^-)$ have a non-empty intersection. Additionally 
the rate of convergence in the definition of the sets $W^\pm$ is controlled by the rate 
function $\phi$.

The origin of the gluing property is the so-called local product structure, introduced by D.V.~Anosov 
for uniformly hyperbolic DS. The local product structure means that for a pair 
of close enough points their stable and unstable manifolds intersect, and the orbit of 
the point of intersection approximates the corresponding semi-trajectories with 
an error decreasing exponentially in time. In our notation this means $\phi(k)=Ce^{-b|k|}$. 
Later in \cite{Bl88} this property has been extended to the global one, but only for 
the uniformly hyperbolic DS.

\begin{remark} (Necessity) The gluing property is necessary for the A,A' and R types of shadowing.
\end{remark}

\proof Already in the simplest case of a single large perturbation the average shadowing 
implies the gluing of any pair semi-trajectories. If $D:={\rm diam}(X)<\infty$ we may choose 
$\phi\equiv D$, otherwise, assuming that the perturbations are bounded by a constant 
$D<\infty$, one shows that during the gluing procedure the gaps between the glued 
segments of true trajectories cannot exceed $KD$, where $K$ depends only on $T$, 
but not on the particular pseudo-trajectory. 
Therefore in the unbounded case we set $\phi\equiv KD$. \qed

Now we are going to check the gluing property for some important classes of 
smooth invertible dynamical systems.

\begin{example}\label{e:hyp} (Affine mapping) 
Let $X:=\IR^d$ with $d\ge1$ with the euclidean metric $\rho$, $A$ be an invertible 
$d \times d$ matrix, and $a\in\IR^d$. Then $Tx:=Ax + a$ is an affine map from $X$ into itself. 
\end{example}

An invertible $d\times d$ real-valued matrix $A$ decomposes the euclidean space $\IR^d$ 
into a direct sum of three linear subspaces $E^s, E^u, E^n$ (stable, unstable and neutral):
$$  E^s:= \{v\in \IR^d:~~ ||A^nv||\le C\la_s^n||v|| ~~\forall n\in\IZ_+\} ,$$
$$  E^u:= \{v\in \IR^d:~~ ||A^{-n}v||\le C\la_u^{-n}||v|| ~~\forall n\in\IZ_+\} ,$$
$$  E^n:= \{v\in \IR^d:~~ ||Av||=||v||\} ,$$
where $0<\la_s<1<\la_u<\infty, C<\infty$. 
In general some of the subspaces $E^s, E^u, E^n$ may be empty. If $E^n=\emptyset$ 
we say that the matrix $A$ is hyperbolic.

If dimension of $E^s$ (or $E^u$) is greater or equal to $d/2$, then it splits $\IR^d$ into 
two $A$-invariant half-hyperplanes ``below'' or ``above'' the corresponding linear subspace. 

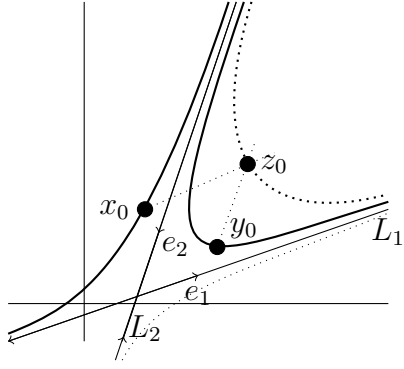
\begin{figure}\begin{center}
\begin{tikzpicture}[scale=0.5]
     \draw [-](0,0) to (10,0); \draw [-](2,-1) to (2,8);
     \draw [-](0,-1) to (10,2.5); \draw [<->](0,-1) to (5,.75); \node at (5,.2){$e_1$}; 
     \draw [-](3,-1) to (6,8);     \draw [>-<](3,-1) to (4,2);    \node at (4.4,1.5){$e_2$}; 
     \node at (10,1.9){$L_1$};  \node at (3.6,-.65){$L_2$}; 
     \draw [thick]  (0,-.8) .. controls (2,0) and (3.5,1) .. (5.8,8);  %Bezier
     \draw [thick]  (6.2,8) .. controls (3.5,0) and (4.5,1) .. (10,2.7);  %Bezier
     \draw [fill] (5.5,1.5) circle (0.2cm);  \node at (6.2,2){$y_0$}; 
     \draw [fill] (3.6,2.5) circle (0.2cm);  \node at (2.8,2.5){$x_0$};
     \draw [dotted](3.6,2.5) to (7,4); \draw [dotted](5.5,1.5) to (6.5,4.3);
     \draw [fill] (6.3,3.7) circle (0.2cm); \node at (7,3.7){$z_0$}; 
     \draw [thick, dotted] (6.4,8) .. controls (5,3.3) and (7,2.1) .. (10,2.9);  %Bezier
     \draw [-](2.82,-1.5) to (6,8);  
     \draw [dotted]  (3,-1.5) .. controls (4,.4) and (6,.75) .. (10,2.4);  %Bezier
\end{tikzpicture}\end{center}
\caption{Hyperbolic affine mapping with $d=2,~E^u=e_1, ~E^s=e_2$. 
             Gluing of $\v{x}$ and $\v{y}$. Two typical trajectories are 
             indicated by thick lines, while the gluing trajectory by a thick dotted line. The lower dotted line 
             corresponds to a trajectory separated from the trajectory $\v{x}$ by both lines $L_1$ 
             and $L_2$, but it can be glued with $\v{x}$ by means of the trajectory $\v{y}$.}
\label{f:hyp-ex} \end{figure}

\begin{proposition} \label{p:hyp}
\begin{enumerate}
\item If $E^n\ne\emptyset$ then $\forall \v{x}~~\exists \v{y}$ with an arbitrary 
         small $\rho(x_0,y_0)\ll1$ which cannot be glued together with a summable 
         accuracy rate $\phi$.
\item If $E^n=E^s=\emptyset$ then $T\in G_s(\phi)$ with $\phi(k):=C\la_u^{-|k|}$.
\item If $E^n=E^u=\emptyset$ then $T\in G_s(\phi)$ with $\phi(k):=C\la_s^{|k|}$.
\item If $E^n=\emptyset$ and $E^u,E^s\ne\emptyset$ then $T \in G_s(\phi)$ with 
        $\phi(k):=C(\la_s^{|k|} + \la_u^{-|k|})$.
\end{enumerate}
\end{proposition}

\proof Let $\v{x}$ be an arbitrary trajectory of the map $T$ and let $0\ne v\in E^n$. 
Consider a trajectory $\v{y}$, defined by the relations 
$y_0:=x_0+v,~~y_k:=T^ky_0~~\forall k\in\IZ$ (see Fig.~\ref{f:hyp-ex}). 
The existence of a trajectory $\v{z}$ which glues $\v{x}, \v{y}$ with 
a summable accuracy rate implies that 
$$ \liminf_{n\to-\infty}\rho(z_n,x_n) + \liminf_{n\to\infty}\rho(z_n,y_n) >0 .$$ 
Therefore $z_0-x_0\in E^u,~~ z_0-y_0=z_0-(x_0+v)\in E^s$ with $v\in E^n$. 
Therefore $(E^u - v)\cap E^s \ne\emptyset$. We came to a contradiction, which 
proves the case (1).

\bigskip

In the case (2) the map $T$ is expanding and the verification of the gluing 
property for a pair of a backward semi-trajectory $\v{x}:=\{\dots,x_{-1},x_0\}$ 
and a forward semi-trajectory $\v{y}:=\{y_0,y_1,\dots\}$ of the map $T$ is straightforward. 
Setting $z_0:=y_0$, we get 
$z_k\equiv z_k~~\forall k\ge0$ and $\rho(z_k,x_k)\le \la_u^{-|k|}~~\forall k<0$. 
Observe, that this is the only possibility to construct the approximating trajectory.

\bigskip

In the case (3) the map $T$ is contracting and hence $T^{-1}$ is expanding. 
Therefore the proof is exactly the same as in the previous case with the only 
change to the inverse time. 

\bigskip

It remains to consider the generic case (4). In this case there is a fixed point 
$O:=(I-A)^{-1}a$. Consider an arbitrary pair of trajectories $\v{x}$ and $\v{y}$ 
(see Fig~\ref{f:hyp-ex}). 
By the assumptions the intersection $I$ of the sets $x_0+E^u$ and $y_0+E^s$ is 
nonempty. Choosing any point $z_0\in I$, for its trajectory we get  
$$ \rho(z_k,x_k)\le C(\la_s^{|k|} + \la_u^{-|k|})~~\forall k<0, $$
$$ \rho(z_k,y_k)\le C(\la_s^{|k|} + \la_u^{-|k|})~~\forall k\ge 0. $$
An important point is that the constant $C=C(A,\rho)$ can be arbitrary large here 
independently on the optimization by the proper choice of the point $z_0\in I$. 
\qed

\begin{remark} If $d=2$, then for the case of a hyperbolic matrix, namely, if 
$\v{x}, \v{y}$ belong to the same $T$-invariant half-hyperplane, then $\v{z}:=\v{y}$ 
glues $\v{x}, \v{y}$ together (see Fig~\ref{f:hyp-ex}). 
Unfortunately, when the dimension of $E^u$ or $E^s$ is greater than 1, 
this simple recipe does not work. 
\end{remark}

\bigskip

\begin{example}\label{e:Anosov} (Anosov diffeomorphism)
Let $X:=\Tor^2$ be a unit 2-dimensional torus and let $T:X\to X$ be a uniformly hyperbolic
diffeomorphism. 
\end{example}

The simplest map that satisfies the above properties is $Tx:=Ax \mod1$, where $A$ is an 
integer matrix with the determinant equal 1 on modulus. For exact definition of the 
uniformly hyperbolic system we refer the reader to numerous publications to the subject 
(see, for example, \cite{An2, Bo, BPSJ, Bl97}).

\begin{proposition} \label{p:Anosov} 
For the map $T:\Tor^2\to\Tor^2$ of the example~\ref{e:Anosov} 
there exists a special (Lyapunov) metric $\rho$ and $\la>1$, for which this system 
satisfies the $G_s(\phi)$ property with $\phi(k):=e^{-\la|n|}~~\forall n\in\IZ$.   
\end{proposition}

This result follows from the global product structure for a hyperbolic system proven in \cite{Bl88}. 
The local version of this property, which asserts the intersection of stable and unstable 
local manifolds of sufficiently close points, is well known (see, for example,  \cite{An2, Bo, BPSJ}). 
In the global version, the locality assumption is dropped. 

It is worth noting that none of these results follow from Proposition~\ref{p:hyp}, 
and not vice versa. Indeed, the diffeomorphism under study is nonlinear and 
the proof of Proposition~\ref{p:Anosov} is based on the construction of 
arbitrary thin Markov partitions and their mixing properties. On the other 
hand, these constructions fail in the example~\ref{e:hyp}. Moreover, 
Proposition~\ref{p:hyp} holds for an arbitrary metric $\rho$ (induced by a norm), 
while Proposition~\ref{p:Anosov} holds only for a special (Lyapunov) metric.

\section{The gluing property for non-uniformly hyperbolic endomorphisms with singularities}\label{s:endo}
Now we are ready to turn to discontinuous and non-invertible mappings. We start 
with a non-invertible version of the example~\ref{e:hyp}.

\begin{example}\label{e:hyp-non} (General affine mapping) 
Let $X:=\IR^d$ with $d\ge1$ with the euclidean metric $\rho$, $A$ be an arbitrary  
$d \times d$ matrix, and $a\in\IR^d$. Then $Tx:=Ax + a$ is an affine non-necessarily  
invertible map from $X$ into itself. 
\end{example}

An general $d\times d$ real-valued matrix $A$ decomposes the euclidean space $\IR^d$ 
into a direct sum of four linear subspaces $E^0,E^s, E^u, E^n$ (kernel, stable, unstable and neutral), 
where 
$$  E^0:= \{v\in \IR^d:~~ Av=0\} ,$$
while three other subspaces were defined in the discussion of the example~\ref{e:hyp}.

\begin{proposition}\label{p:hyp-non} The map $T\in \cS(A,A')$ if and only if 
$E^n=\emptyset$ and either $E^s=\emptyset$ or $E^u=\emptyset$.
\end{proposition}

\proof This claim differs from Proposition~\ref{p:hyp} in two ways. 
First, we need to take into account the presence of the subspace $E^0$, and second, 
we need to prove necessary part.

The first part can be proven as follows. Observe that the restriction of the map $T$ to 
the invariant linear subspace $E^s + E^u + E^n$ satisfies Proposition~\ref{p:hyp}. 
By definition, $TE^0=0$. Therefore when checking the gluing property it is enough 
to repeat the construction of $\v{z}$ from the proof of Proposition~\ref{p:hyp} 
for the points $x_0,y_1$ instead of $x_0,y_0$. 

To prove necessary part, observe that the $T$ image of each perturbation in the 
direction $E^0$ vanishes. Therefore under the assumption of small on average 
perturbations their contribution to the approximation error is small as well. \qed

\bigskip

So far the maps in all the examples under consideration were linear, 
and when the gluing property was satisfied for them, the rate function was exponential. 
The following example demonstrates that this is not necessary.

\begin{example}\label{ex:non-unif} (Nonuniform hyperbolicity) $X:=[0,1],~\alpha,\beta\ge0$, $0<c<1$
$$ Tx:=\function{x(1+ax^\alpha) &\mbox{if } x\le c \\ 
                       1 -  (1-x)(1+b(1-x)^\beta) &\mbox{if } x>c} .$$ %$0<c,~0\le\alpha\le1$. 
\end{example} 

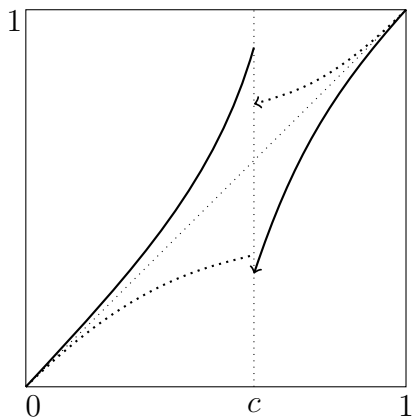
\begin{figure}\begin{center}
\begin{tikzpicture}[scale=0.5]
     \draw [-](0,0) to (10,0) to (10,10) to (0,10) to (0,0);  
     \node at (0.2,-0.5){$0$}; \node at (10,-0.5){1}; \node at (-0.3,9.7){1}; 
     \draw [dotted] (0,0) to (10,10); \draw [dotted] (6,0) to (6,10); \node at (6,-0.5){$c$};
     \draw [thick,black]  (0,0) .. controls (3,3.25) and (5,5.5) .. (6,9);  %Bezier
     \draw [thick,black,->]  (10,10) .. controls (8,7.75) and (7,6) .. (6,3);  %Bezier
     \draw [thick,dotted]  (0,0) .. controls (3,3) and (5,3.2) .. (6,3.5);  %Bezier
     \draw [thick,dotted,->]  (10,10) .. controls (8,8) and (7,7.8) .. (6,7.5);  %Bezier
\end{tikzpicture}\end{center}
\caption{Graph of the nonuniform hyperbolic map $T$ with $\alpha,\beta>0$. 
Thick lines correspond to $a,b>0$, while thick dotted lines to $a,b<0$.}\label{f:neut} \end{figure}

This example is of particular interest because it exhibits very different behaviors for 
different regions in the parameter space $(a,b,c,\alpha,\beta)$.  
Observe, that if $\alpha,\beta=0,~a,b>0$ we are getting a piecewise expanding map, 
while for $\alpha,\beta=0,~a,b<0$ this is a contracting map. The most interesting 
situation (nonuniform hyperbolicity)) corresponds to the case $\alpha,\beta>0$, 
when there are two neutrally expanding fixed points at the end-points $0$ and $1$. 
From the point of view of ergodic theory, there are additional peculiarities in the 
properties of invariant measures. Assume that the equality (\ref{e:ends}) below holds true. 
Then for $0<\alpha,\beta<1$ there is an absolutely continuous probabilistic invariant measure, 
while for $\alpha,\beta>1$ there are exactly two ergodic probabilistic invariant measures -- 
Dirac-measures at points 0 and 1 (apart from $\sigma$-finite absolutely continuous 
probabilistic invariant measures). 
As we will see, this important distinction is reflected by the presence or absence 
of the shadowing property.

An important feature for our analysis of shadowing, common for all admissible choices 
of the parameters $(a,b,c,\alpha,\beta)$, is the presence of at least two periodic 
trajectories -- fixed points at 0 and 1. 

\begin{proposition}\label{p:non-unif} The map $T$ in the example~\ref{ex:non-unif} 
(see Fig.~\ref{f:neut}) satisfies 
\begin{enumerate}
\item if $\alpha,\beta\ge0$ and $T\in G_{s/w}(\phi)$ then %
   \beq{e:ends}{c(1+ac^\alpha)=(1-c)(1+b)^\beta=1 ;} %
\item if $\alpha=\beta=0$ and (\ref{e:ends}) holds true, then $T\in G_s(\phi)$ 
   with an exponentially decreasing rate function $\phi$;
\item if $0<\alpha,\beta<1$ and (\ref{e:ends}) holds true, then $T\in G_w(\phi)$
        with the summable rate function 
        $\phi(k):=\function{C|k|^{-\gamma} &\mbox{if~} k\le0 \\
                                                    0 &\mbox{if~} k\ge0 }$, 
        where $C=C(a,b,c,\alpha,\beta)<\infty,~\gamma>1/\min(\alpha,\beta)$; 
\item if $\alpha,\beta>1$ and $ab\ne0$, then the strong gluing property (\ref{e:glu}) 
         with a summable rate function $\phi$ cannot hold.
\end{enumerate}
\end{proposition}

\begin{remark}
(a) $\alpha=\beta=0$ means that the function $T$ is piecewise linear.\\
(b) Equality (\ref{e:ends}) is equivalent to $T[0,c]=[0,1],~T(c,1]=(0,1]$, which implies that $a,b>0$.
\end{remark}

\proof If any version of the gluing property holds, then any pair of trajectories of the map $T$ 
may be glued together. Consider a pair of trajectories $\v{x}:=\{0\},~ \v{y}:=\{1\}$. 
Here by $\{0\}$ and $\{1\}$ we mean trajectories staying at fixed points $0$ and $1$ 
correspondingly $\forall t\in\IZ$. 

There are several possibilities:
\begin{enumerate}
\item $a,b<0$ or $a,b>0$. Then for $Tc<1$ then the forward trajectory, starting from a point 
$u$ belonging to a neighborhood of $0$, cannot assume a value greater $Tc<1$ (here we are 
using monotonicity of the branches of $T$). 
Similarly, if $T(c-0)>0$ any backward trajectory, starting from a point $v$ belonging 
to a neighborhood of $0$, cannot assume a value smaller $T(c-0)>0$. 
Therefore no gluing of $\v{x}$ and $\v{y}$ with a summable rate function may take place.
\item $a<0<b$. Then a forward trajectory, starting from a point $u$ belonging to a 
neighborhood of $0$, cannot get out of this neighborhood. Therefore we came to the same conclusion. 
\item $a>0>b$. Similarly to the previous item, but one needs to consider the neighborhood of $1$.
\end{enumerate}

This analysis proves item (1).

\bigskip

To prove item(2) consider two arbitrary trajectories $\v{x}, \v{y}$ of this map. 
We define the gluing trajectory $\v{z}$ as follows: 
$z_0:=y_0, z_k:=T^ky_0$ for $k>0$, and $z_{k-1}:=T^{-1}_{x_{k-1}}z_k$ for $k\le0$. 
In fact, since the map $T$ is piecewise expanding, there are no other options for $z_0$, 
similarly to the example~\ref{e:hyp}.

By this construction $z_k=y_k ~~\forall k\ge0$, while for negative $k$ the distances 
between $z_k$ and $x_k$ decrease at exponential rate, since each time we are applying for their 
calculations the same inverse branch of the expanding map $T$. 

Nevertheless, observe that the distance between $Tx_0$ and $Ty_0$ (i.e. the gap between 
the backward trajectory $\v{x}$ and the forward trajectory $\v{y}$) might be arbitrary 
close to $1$ when the points $x_0$ and $y_0$ are close to the point $c$. 

Thus if $\min(a,b)>0$ we get the summable rate function 
$$ \phi(k):=\function{(1+\min(a,b))^{k} &\mbox{if~} k\le0 \\
                                                  0 &\mbox{if~} k\ge0 } .$$   
Item (2) is proven.

\bigskip

%\bcr{This result demonstrates that it is possible that $T\in (A+A)$ while $T\not\in (U+A)$. [??]}

In all situations considered so far the rate function $\phi$ has exponential tails. In general this is 
absolutely not the case, which will be demonstrated in the case $\alpha,\beta>0$. 
Moreover, we will show that the neutrally expanding map, considered in this example, 
for $0<\alpha,\beta<1$ satisfies only the weak gluing property (\ref{e:glu-w}), but 
the strong one (\ref{e:glu}) breaks down, which excludes the uniform (U+U) shadowing property.

If $\alpha,\beta>0$ the map $T$ is uniformly piecewise expanding everywhere except the 
neighborhoods of the neutral fixed points $0$ and $1$. 
Therefore the analysis of the forward part of the gluing property does not differ much 
from the situation in item (2), namely, we choose $z_0:=y_0$. 
Still we need to show that proper chosen pre-images of the point $z_0$ will approximate 
the backward semi-trajectory well enough. 

To this end we need some estimates of the rate of convergence in backward time of a map 
with a neutral fixed point obtained in \cite{Bl22}.  

\begin{lemma}\label{l:neut} \cite{Bl22} Let $\tau(v):=v+Rv^{1+\alpha},~R>0, \alpha\ge0,~v\ge0$. Then \\
(a) $\tau^{-n}(v) \le K n^{-\gamma}~~\forall v\in[0,1],~n\in\IZ_+$ 
     and some $K<\infty,~\gamma>1/\alpha$.\\
(b) $\tau^{-n}(v) \ge Kv n^{-\gamma}~~\forall v\in[0,1],~n\in\IZ_+$ 
     and some $K<\infty,~\gamma<1/\alpha$.\\
(c)  if $\alpha=0$ then $\tau^{-n}(v) \le (1+R)^{-n}v~~\forall v\in[0,1],~n\in\IZ_+$.\\
\end{lemma}

Applying the assertion (a) of Lemma~\ref{l:neut} to the inverse branches of the map $T$, 
and using that $0\le u\le1$ we get %
\beq{e:neut}{ \rho(T^{-n}u, \{0,1\}) < C n^{-\gamma}~~\forall u\in X,~n\in\IZ_+ ,} %
where $C=C(a,b,c,\alpha,\beta)<\infty$, $\gamma>1/\min(\alpha,\beta)$ and 
$$ \rho(u,A):=\inf_{a\in A}\rho(u,a), ~~\rho(u,v):=|u-v| .$$ 
Now we are ready to estimate $\rho(x_{-n},z_{-n})$ for $n\in \IZ_+$. 
By the triangle inequality, using (\ref{e:neut}), we get 
$$ \rho(x_{-n},z_{-n}) \le \rho(x_{-n},0) + \rho(0,z_{-n}) 
     \le 2C n^{-\gamma} .$$
Therefore, since $z_n\equiv y_n~\forall n\in\IZ_+$, 
$\phi(k):=\function{2C|k|^{-\gamma} &\mbox{if~} k\le0 \\
                                              0 &\mbox{if~} k\ge0 }$ 
defines the rate function for the (weak) gluing property (\ref{e:glu-w}).  
Moreover, if $\alpha<1$ then $\gamma>1$ and hence $\phi$ is summable. 
This proves item (3).

\bigskip

Finally, the fact that $T\not\in G_s(\phi)$ for any summable $\phi$  
follows from Lemma~\ref{l:neut}(b), which proves item (4). \qed

%\begin{example}\label{ex:2}\bcr{(Kan-Yorke)}
%Let $X:=[0,1], Tx:=\function{ax &\mbox{if } x<c \\  b(1-x)  &\mbox{otherwise}}$. 
%\end{example}

\section{The gluing property for symbolic dynamics}\label{s:symb}
\begin{example}\label{ex:4} (Symbolic dynamics)
Let $T:X\to X$ be a map and let $\{X_i\}_{i=1}^r$ be a partition of $X$.\footnote{The partition 
    needs not to be finite, i.e. it is possible that $r=\infty$.} 
Then each (pseudo-)trajectory $\v{x}$ of the map $T$ may be coded by a bi-infinite sequence 
$\v{s}$ of elements from the alphabet $\cA:=\{a_1,a_2,\dots,a_r\}$ according to the 
rule: $s_i:=a_k$ if $x_i\in X_k$. This gives a symbolic description of the dynamics, 
governed by a binary {\em transition} matrix $\pi$, where 
$\pi_{ij}=1$ iff $TX_i\cap X_j\ne\emptyset$. The set of admissible sequences 
(corresponding the transition matrix $\pi$) we denote by $\v{\Sigma}_\pi$, while 
the set of all sequences from the alphabet $\cA$ by $\v{\Sigma}:=\cA^\IZ$. When coding a 
pseudo-trajectory $\v{y}$ of the map $T$ we get a sequence $\v{s}\notin \v{\Sigma}_\pi$. 
The question is, is there is an admissible sequence that ``approximates'' $\v{s}$?
\end{example}
This example differs from the previous ones in two important points. First, the ``amplitude''  
of the perturbation takes only a finite number of values and thus cannot be uniformly small. 
Second, despite the transition between letters from the alphabet $\cA$ can be described 
in terms of a map, but this map is multi-valued and so our results about shadowing cannot 
be applied in this setting directly. 

To this end we consider a shift map $\sigma:\v{\Sigma}_\pi\to \v{\Sigma}_\pi$ with somewhat 
unusual perturbations. For a pair of admissible sequences $\v{s},\v{u}$ by the perturbed 
sequence we mean $\v{w}\in \v{\Sigma}$ such that 
$w_i:=\function{s_i &\mbox{if } i<0 \\  u_i &\mbox{otherwise}.}$ 
In words, we preserve the elements with negative indices of $\v{s}$, but change the ones 
with nonnegative indices to $\v{u}$. Obviously, this is what happens under the perturbations 
in the Example~\ref{ex:4}. In the space of all sequences $\v{S}$ we consider a metric 
$\rho(\v{s},\v{u}):=\sum_{k=-\infty}^\infty 2^{-|k|} 1_{s_k}(u_k)$, 
where $1_a(b):=\function{0 &\mbox{if } a=b \\  1 &\mbox{otherwise}.}$. 

Clearly, the average shadowing in symbolic dynamics is equivalent to the average shadowing 
for the shift map $\sigma$ acting in the metric space $(\v{\Sigma},\rho)$.

\begin{proposition}\label{p:symb} Symbolic dynamics belongs to the class $G_s(\phi)$ with 
a finitely supported (and hence summable) function $\phi$ if and only if 
$\exists M\in\IZ_+$ such that the $M$-power of the transition matrix (i.e. $\pi^M$) is positive. 
Therefore under the latter assumption the symbolic dynamics belongs to the class $\cS(R,A)$.
\end{proposition}

Before to give the proof, let us demonstrate examples of the transition matrices, 
leading to the average shadowing 
$\pi:=\left(\begin{array}{cc} 1   & 1 \\  1 & 0 \end{array} \right)$, or its absence 
$\pi:=\left(\begin{array}{cc} 0   & 1 \\  1 & 0 \end{array} \right)$.

\proof We start with the case of the finite alphabet $\cA$, i.e. $r<\infty$. 
According to the assumption $\pi^M>0$ during time $M$ we can go from any element of 
the alphabet $\cA$ to any other. So setting the function $\phi$ equal to the indicator function 
of the integer segment $[-M,M]$, we obtain the desired result. Namely, this demonstrates 
that each pair of admissible trajectories may be glued together with accuracy $\phi$. 

\bigskip

To verify the necessary part, assume that the condition of the theorem is not satisfied
and $\not\exists M:~~\pi^M:=(\pi_{ij}^{(n)})>0$. It follows that
$$ \forall n\in\IZ_+ ~~\exists i_n,j_n:~~\pi_{i_n,j_n}^{(n)}=0 .$$
Indeed, if this were not the case, then there would be $k\in\IZ_+$ such that $\pi^k>0$.
But in this case $\pi^n>0~\forall n>k$, which contradicts the assumption.

As we have already noted, the transition matrix $\pi$ induces a multivalued mapping
of the alphabet into itself by the following formula: $\pi a:= \{b\in\cA:~~\pi_{ab}>0\}$.
In these terms, since $\pi$ is not strictly positive and $\cA$ is finite,
it follows that for some $N\in\IZ_+$ there is a partition $\cA:=\sqcup_i\cA_i$
into non-empty $\pi^N$-invariant subsets, i.e. $(\pi^N)^{-1}\cA_i=\cA_i~\forall i$.
It follows from this that for $x_0,y_0$ belonging to different elements of this
partitions corresponding to admissible sequences $\v{x},\v{y}$ do not
intersect, i.e. $x_i\ne y_i~\forall i\in\IZ$. Therefore, there is no ``gluing'' of them
with any summable function $\phi$.

\bigskip 

It remains to prove that the claim of theorem remains valid for $r=\infty$. 
The sufficient part follows from the same argument as in the case $r<\infty$. 
To prove the necessary part, consider a ``circumcised'' finite alphabet 
$\t\cA:=\{a_1,a_2,\dots,a_{\ell-1},\t{a}_{\ell}\}$ with $\ell<\infty$ letters, such that 
$\{a_1,a_2,\dots,{a}_{\ell-1}\}\subset\cA$, while the letter $\t{a}_{\ell}$ corresponds 
to the union of the remaining elements of the partition $\cup_{i\ge \ell}X_i$. 
Then considering the corresponding ``circumcised'' transition matrix $\t\pi$, we 
get the shift map $\t\sigma: \v{\Sigma}_{\t\cA}\to \v{\Sigma}_{\t\cA}$, for which all results 
obtained in the first part of the proof are applicable. Therefore the average shadowing 
for $\t\sigma$ holds iff $\exists \t{M} = \t{M}(\ell) <\infty:~\t\pi^{\t{M}}>0$. 

According to the validity of the average shadowing for $\t\sigma$, the  inequality 
$\pi^M>0$ does not hold $\forall M\in\IZ_+$ for the complete infinite alphabet $\cA$ 
if and only if $\t{M}(\ell)\toas{\ell\to\infty}\infty$. If the latter happens, one 
constructs a pair of admissible trajectories, having elements from $\cA$ with arbitrary 
large indices, which cannot be glued together with the summable accuracy rate $\phi$. 
We came to the contradiction. 
\qed

%%%%%%%%%%%%%%%%%%%%%%%%%%%

%\small
%\newpage

\end{document}